\documentclass[12pt]{article}
\usepackage{amsmath, amsthm, amssymb}
\usepackage{latexsym,amsfonts,amssymb,amsthm,amsmath, amscd,pifont,makecell}
\usepackage[titles]{tocloft} 
\usepackage[square, sort, numbers, authoryear]{natbib}
\usepackage[T1]{fontenc}
\usepackage{babel}

\bibpunct{(}{)}{;}{a}{,}{,}
\usepackage{epsf,amsfonts,amsmath}
\usepackage{hyperref} 
\renewcommand{\thesection}{\arabic{section}.}

\newtheorem{thm}{Theorem}[{section}]
\newtheorem{lemma}{Lemma}[section]

\newtheorem{defn}{Definition}[{section}]

\def\th{\theta}
\def\Th{\Theta}

\newcommand{\be}{\begin{equation}}
\newcommand{\enq}{\end{equation}}

\begin{document}

\begin{titlepage}
\title{\bf The Equivariance Criterion in a Linear Model for Fixed-X Cases}
\author{Daowei Wang\\
Paula and Gregory Chow Institute for Studies in Economics\\
Xiamen University\\
Xiamen, Fujian China PRC\\
ORCiD: 0009-0007-3428-723X\\
   \texttt{wangdw97@stu.xmu.edu.cn}
   \and
   Mian Wu\\
   Department of Computer Science and Engineering\\
   The Ohio State University\\
   Columbus, OH, USA \\
   ORCiD:  0009-0001-2116-1773\\
   \texttt{wu.6666@osu.edu}
   \and
   Dr. Haojin Zhou\thanks{Corresponding author: haojin\_zhou@hotmail.com} \\
Academy of Pharmacy \\
Xi’an Jiaotong - Liverpool University (XJTLU) \\
Suzhou, Jiangsu China PRC\\
   ORCiD: 0000-0003-0802-099X\\
      \texttt{haojin.zhou@xjtlu.edu.cn}
}
\date{}
\maketitle
\begin{abstract}
The field of machine have seen rising applications of equivariance criterion. However, there is no systematic way to justify its usage, including why it
works, whether there is an optimal solution and if so, what form it carries. In this article, we explored the usage of equivariance criterion in a normal
linear model with fixed-$X$ and extended the model to allow multiple populations, which, in turn, leads to a multivariate invariant location-scale
transformation group, compared than the commonly used univariate one. The minimum risk equivariant estimators of the coefficient vector and the diagonal
covariance matrix were derived, which were consistent with literature works. This work serves as an early exploration of the usage of equivariance
criterion in machine learning, where we confirmed that the least square approach widely used in machine learning indeed carries optimality in some sense at
least in the framework of estimation.

Meanwhile, the problems can be shown to be equivalent to a mixture from $p$ independent normal samples and via the principle of functional equivariance, an alternative proof can be derived. However, such an approach carries its own limitation with a strong tie to equivariance criterion.

\end{abstract}


\maketitle

\end{titlepage}

\section{Introduction}
\renewcommand{\thesection}{\arabic{section}}

Consider a general linear model with a response $Y \in \mathbb{R}$ and a covariate vector ${\bf X} \in \mathbb{R}^p$. The aims of the linear model are to
build a linear functional relationship between $Y$ and ${\bf X}$
from paired observations, $({\bf X}_1,Y_1),\dots,({\bf X}_n,Y_n)$, and predict a future response, $Y_0$, given a covariate, ${\bf X}_0$. Predominately in
literature, a fixed-$X$ assumption has been used with the following
understandings:

(i) The covariate values are fixed before sampling while the only randomness comes from the responses $Y_1,\dots,Y_n$;

(ii) There is no other relationship among those covariate values.

Meanwhile, random-$X$ cases have been investigated, e.g., (\citet{BreimanSpector1992}), where the assumption is to set the covariate values sampled from a
random vector.

In this article, we will start with a fixed-$X$ as a combination of $p$ linearly independent rows each repeated $n_i\geq 1$ times, i.e.,
\[
X = (\bold{X}_1,\ldots,\bold{X}_1, \bold{X}_2,\ldots,\bold{X}_2, \dots, \bold{X}_p,\ldots,\bold{X}_p)',
\]
where $\sum n_i=n$. One can argue that this is the general form for a fixed-$X$ case as the rank of the design matrix $X$ is $p$ allowing only $p$ distinct
rows. The difference between a fixed-$X$ and a random-$X$ lays on
the fact whether $X$ is fully known before sampling. Naturally, one should set the number of parameters as the rank of $X$, $p$. Otherwise if there are more
parameters than $p$, then the model carries redundant parameters, some of which may not be identifiable and thus not estimable. On the other hand, if there
are fewer parameters than $p$, then more parameters will be needed or some of $X$ may be considered as unknown before sampling and hence $X$ is random and
out of scope of this article. To derive the optimal solution, we will introduce equivariance criterion, which has been widely considered in the linear model,
e.g., \citet{Rao1965,Eaton1989} and the distinctions between fixed-$X$ and
random-$X$ cases are clearly noticed (\citet{Rao1973}, \citet{RossetTibshirani2020}).

As a principle of symmetry, equivariance criterion is proposed in literature as a logic constraint on the solutions to derive the optimal one. Compared to
the unbiasedness, it is more focused on the symmetry of the problem
and self-consistency of the solutions. In the era of machine learning, it has received rising attention due to widely existing symmetry in the applications.
\citet[Chap. 3]{LehmannCasella1998} has given a detailed
discussion of equivariance criterion in the location-scale family while \citet[Chap. 6]{Berger1985} has presented the theory in a decision theoretic
framework. Besides those two classical textbooks, other important
references are as follows: \citet{HoraBuehler1966, Eaton1989, Wijsman1990}.

Equivariance criterion consists of two principles: (i) functional equivariance, which states the action in a decision problem should be consistence across
different measures and (ii) formal equivariance, which requires the
decision rule to be of the same form for two problems with the identical structure.

Formally, a decision problem is described by $({\cal X},{\cal P},\Th,{\cal D},L)$, where ${\cal X}$ is a sample space, $\cal P$ is a family of distributions
with $\th$ as the parameter or the true state of nature, $\Th$ is
the parameter space, $\cal D$ is a decision space and $L$ is a loss function on ${\cal D} \otimes \Th$.

Without loss of generalities, one starts with the identifiability on both the distribution family and loss function, which, for linear model, suggests that
there is no redundant $\boldsymbol\beta\in \mathbb{R}^p $ and
$d\in \mathbb{R}^p$. The principle of functional equivariance criterion requires preservance of the model and invariance of loss function under a group of
1-1 and onto transformations, which, are defined as follows.

\begin{defn}
(Preservance of the Model) The distribution family ${\cal F}$ is said to be invariant under $G$ if for each $g\in G$, and $\theta \in \Theta$, there exists
$\theta' \in \Theta$ such that
\be
X\sim f(x|\th)\Rightarrow g(X) = X' \sim f(x'|\th')
\label{Preservance}
\enq
\end{defn}

\begin{defn}
(Invariance of Loss Function) The loss function $L(d, \th)$ is invariant under $G$ if for each $g\in G$ and each $d\in {\cal D}$ there exists
$d_*=\tilde{g}(d) \in {\cal D}$ such that
\[
L(d, \th) = L(\tilde{g}(d), \bar g(\th)) \mbox{ for all } \th \in \Th.
\]
\end{defn}

\begin{defn}
A decision problem $({\cal X},{\cal F},\Th, {\cal D}, L)$ is invariant under a group of transformations, $G$, if
$G$ preserves $\cal{P}$ and the loss function is invariant under $G$.
\end{defn}
Under the structure of an invariant decision problem, we can induce equivariance criterion.

\begin{defn}
A decision rule $\delta(X)$ is said to be equivariant under $G$ if
\be
\delta(g(x))=\tilde{g}(\delta(x)) \quad \hbox{for all } g \in G \hbox{ and } x \in {\cal X}.
\label{equiv}
\enq
\end{defn}

In this way, equivariance criterion imposes a constraint on the possible decision rules that one can use reasonably to derive the optimal decision rule. In
this article, we will pursue the best decision rule (minimum risk equivariant rule, MRE rule) in the sense of minimizing the risk function as the expectation
of the loss function over ${\bf X}$, which is a function of the parameter and the decision rule.

\citet{LehmannCasella1998} has applied equivariance criterion implicitly in linear model for fixed-$X$ cases, where the most important feature is to
transform the problem to a canonical form and derive the best equivariant
estimators for the coefficient vector and common variance. \citet{WuYang2002} has discussed the existence of best equivariant estimators for the coefficient
vector and the covariance matrix (in three forms) in normal
linear models with fixed-$X$s and derived the forms when they exist. \citet{KurataMatsuura2016} has derived the best equivariant estimator of regression
coefficients in a seemingly unrelated regression (SUR) model with a
known correlation matrix. Further, \citet{MatsuuraKurata2020} derived the best equivariant estimator of the variance matrix in a SUR model. It is noted that
most literature works assume a single population for the linear
model and thus apply a common invariant transformation on the responses to derive the optimal equivariant estimators. However, as in experimental design, one
usually views each distinct covariate vector as an independent
population from any other and a SUR model uses a natural multivariate extension of the common linear model. Therefore, a multivariate extension of the common
linear model combined with equivariance criterion warrants
further investigation, especially from the fundamental logic of equivariance criterion.

This article focuses on applying the logic of equivariance criterion to the fixed-$X$ linear model, where the linear model is extended to allow multiple
populations instead of a single one and the invariant transformation
is distinct for each population. In Section 2, we derive the best equivariant estimators for the coefficient vector and the diagonal covariance matrix in a
normal linear model specially tuned for equivariance criterion. In Section 3, we explore alternative proofs via the principle of functional equivariance and discuss such an approach.
Section 4 is devoted to some concluding remarks \& future work.

\setcounter{equation}{0}

\renewcommand{\thesection}{\arabic{section}.}
\section{Equivariance in the Linear Model}
\renewcommand{\thesection}{\arabic{section}}
\setcounter{equation}{0}
Consider the linear regression model ${\bf y}=X {\bf \boldsymbol\beta}+{\bf \epsilon}$, where ${\bf y}$ is an $n\times 1$ vector, $X$ is an $n\times p$
matrix and ${\bf \epsilon}$ is the noise vector. To derive the best
equivariant estimators in the model, we will walk through the basic elements of equivariance criterion first and then provide a linear model tuned for
equivariance criterion starting from the basic concept of a population.

{\bf Preservance of the Model}: Without loss of generalities, we assume that the design matrix $X$ is predetermined and of full rank with $n\geq p+1$, where
only the response ${\bf y}$ is random and thus transformations
only on ${\bf y}$ will be considered. For a fixed-$X$ linear model, it can be viewed that the samples actually come from different populations as in the
experimental design, which imposes a restriction over the choice of
the transformation group, in addition to the specification of the model. Equivariance criterion implicitly requires the transformations over the samples from
the same population to be identical while being distinct for
samples from different populations. Therefore, it is of essence to determine the number of populations inside a linear model.

For the fixed-$X$ cases, one can argue that there can only be $p$ populations inside the linear model as the rank of the design matrix is $p$. In this
regard, we will assume that $\epsilon_{i_j}$'s are independently
distributed as a normal distribution with mean $0$ and unknown variance $\sigma_i^2$, $i=1,\ldots,p$. Thus, we have the normal linear model tuned for
equivariance criterion as follows,
\be
\begin{split}
{\bf y}&=\left( \bold{y}_1',\bold{y}_2',\cdots,\bold{y}_p'\right)'\\ 
\text{ with X }&=({\bf X_1},\ldots,{\bf X_1},{\bf X_2},\ldots,{\bf X_2},\ldots,{\bf X_p},\ldots,{\bf X_p})',\\
{\bf \boldsymbol\beta}&=(\beta_1,\beta_2,\ldots,\beta_p)',\\
\text{ and } \Sigma &=\text{diag}(\sigma_1^2,\ldots,\sigma_1^2,\sigma_2^2,\ldots,\sigma_2^2,\ldots,\sigma_p^2,\ldots,\sigma_p^2).
\end{split}
\label{NormalModel}
\enq
Note that here, it is natural to assume independent variances for each population in contrast to the traditional linear model where all the variances for
${\bf y}$ are assumed to be identical. In this case, we will have
the usual identifiability of the model with respect to ${\bf \boldsymbol\beta}$ and $\Sigma$. Since
$\Sigma=\text{diag}(\sigma_1^2,\ldots,\sigma_1^2,\sigma_2^2,\ldots,\sigma_2^2,\ldots,\sigma_p^2,\ldots,\sigma_p^2)$ contains only
$p$ parameters, one can simplify the problem to consider $\Sigma_p=\text{diag}(\sigma_1^2,\ldots,\sigma_p^2)$.

The transformation group $G$ keeping the model invariant is of the form {$g({\bf y})=C{\bf y}+{\bf a}$, with
$C=\text{diag}(c_1,\ldots,c_1,c_2,\ldots,c_2,\ldots,c_p,\ldots,c_p)$, and ${\bf
a}=(a_1,\ldots,a_1,a_2,\ldots,a_2,\ldots,a_p,\ldots,a_p)'$, $c_i>0$, $a_i\in\mathbb{R}, i=1,\ldots,p$. (Both $c_i$ and $a_i$ are repeated $n_i$ times.) Here
we use independent transformations for those populations, compared to
the literature where most are using one single common transformation, $g_c({\bf y})=c{\bf y}+a,c>0,a\in\mathbb{R}$. To facilitate discussion, we denote that
$X_p=({\bf X_1,X_2},\ldots,{\bf X_p})'$,
 $C_p=\text{diag}(c_1,\ldots,c_p)$, ${\bf a}_p=(a_1,\ldots,a_p)'$ and thus $X=KX_p,{\bf a}=K{\bf a}_p$,
with
$$
K=\left[ \begin{matrix}
	\bold{1}_{n_1}&	\bold{0}			&\cdots	&\bold{0}	\\
	\bold{0}& \bold{1}_{n_2}	&\cdots	& \bold{0}\\
	\vdots&\vdots			& \ddots   & \vdots  &		\\
	\bold{0}	&	\bold{0}	& \cdots &\bold{1}_{n_p}\\
	\end{matrix} \right] _{n\times p}.
$$
The corresponding transformation group $\bar{G}$ on the parameter space is of the form \[{\bar{g}({\bf \boldsymbol\beta},\Sigma_p)=(X_p^{-1} C_p X_p
\boldsymbol\beta+X_p^{-1} {\bf a}_p, C_p \Sigma_p C_p')}.\] It can
be shown that the parameter space is transitive under $\bar{G}$ as $n \geq p+1$.

In the linear model, the usual targets of interest are the coefficient vector $\boldsymbol\beta$ and the covariance matrix $\Sigma$. We will discuss those
two separately in the context of {\bf Invariance of the
Loss Function}.

To derive an MRE decision rule, \citet[Chap. 3]{LehmannCasella1998} have used maximal invariants to characterize all the equivariant estimators and then
minimized the constant risk when $\bar{G}$ is transitive for
location-scale families. We will follow such an approach to derive the best estimators as follows.

\subsection{\texorpdfstring{Estimation of the Coefficient Vector \textbf{$\beta$}}{Estimation of the Coefficient Vector beta}}

 \citet{Staudte1971, ZhouNayak2014} have introduced a method to construct an invariant loss function based on the target of interest. Following their method,
 one can build an invariant loss function as follows, $L_{\bf \boldsymbol\beta}({\bf d,\boldsymbol\beta})=({\bf d}-{\bf \boldsymbol\beta})^T X_p^T
 \Sigma_p^{-1} X_p ({\bf d}-{\bf \boldsymbol\beta})$. This is an extension of
the one implicitly used in the equivariance literature (Theorem 4.3 (b) in \citet[Chap. 3]{LehmannCasella1998} and (1.4) in \citet{WuYang2002}). Thus one
will have the corresponding invariant transformation on the decision space as $\tilde{g}({\bf d})=X_p^{-1} C_p X_p {\bf d} + X_p^{-1} {\bf a}_p$ and
equivariant criterion as ${\bf \delta}(g({\bf y}))={\bf \delta}(C{\bf y}+{\bf
a})=\tilde{g}({\bf \delta}({\bf y})) =X_p^{-1} C_p X_p {\bf \delta}({\bf y})+ X_p^{-1} {\bf a}_p$. It can be shown that those invariant transformations form
a group $\tilde{G}_{\bf \boldsymbol\beta}$ and the least square estimator
is essentially the vector of sample means as in the following Lemma \ref{LSM} and thus equivariant under the group $\tilde{G}_{\bf \boldsymbol\beta}$.
\begin{lemma}
\label{LSM}
For the fixed-$X$ linear model in \ref{NormalModel}, the least square estimator, $(X'X)^{-1}X'{\bf y}=X_p^{-1}{\bf \bar{y}}$, with ${\bf
\bar{y}}=(\bar{Y}_1,\ldots,\bar{Y}_j,\ldots,\bar{Y}_p)'$, where $\bar{Y}_i=\frac{1}{n_i}\sum_{j=N_{i-1}+1}^{N_i}Y_j$ is the sample mean for each population with $N_0=0$, $N_i=\sum_{j=1}^{i}n_j$ for $i=1,\ldots,p$.
\end{lemma}
\begin{proof}
Since $K'K=\text{diag}(n_1,\ldots,n_j,\ldots,n_p)$ and $K'{\bf y}=(n_1\bar{Y}_1,\ldots,n_j\bar{Y}_j,\ldots,n_p\bar{Y}_p)'$, the least square estimator satisfies
that
\begin{equation*}
	(X'X)^{-1}X'{\bf y}=(X_p'K'KX_p)^{-1}X_p'K'{\bf y}=X_p^{-1}(K'K)^{-1}K'{\bf y}=X_p^{-1}{\bf \bar{y}}.
\end{equation*}
\end{proof}

{\bf A Characterization of Equivariant Estimators.} It follows from \citet{LehmannCasella1998} that a characterization of equivariant estimators can be given
by
$$
\delta({\bf y})=(X'X)^{-1}X' {\bf y}+
X_p^{-1}S({\bf y})\omega(\bf{z}),
$$
where $S({\bf y})=\text{diag}(s_1,\ldots,s_p)$, $s_i=\sqrt{\frac{1}{n_i-1}\sum_{j=N_{i-1}+1}^{N_i} (Y_j-\bar{Y}_i)^2}$ is
the sample deviation for each population,
${\bf
z}=(\bold{z}_1',\ldots,\bold{z}_p')'$, with $\bold{z}_i=((Y_{N_{i-1}+2}-Y_{N_i})/(Y_{N_{i-1}+1}-Y_{N_i}
),\ldots,(Y_{N_i-1}-Y_{N_i})/(Y_{N_{i-1}+1}-Y_{N_i}),\text{sign}(Y_{N_{i-1}+1}-Y_{N_i}))'$, is a  maximal invariant, and $\omega(\bold{z})$ is a $p$-dimension vector.

{\bf Best Equivariant Estimator.} To derive the best equivariant estimator minimizing the risk function, first one will use the fact that the risk function
is constant for any equivariant estimator as the parameter space is transitive. Then one can show that $(X'X)^{-1}X' {\bf y}$ is independent of ${\bf z}$.
Meanwhile $(X'X)^{-1}X' {\bf y}$ depends on $\bold{y}$ only via $\bar{\bold{y}}$, which is independent of  $S({\bf y})$.
With above arguments, one can show the main result as follows.
\begin{thm}
\label{MRECoeff}
For the fixed-$X$ linear model in \ref{NormalModel}, the least square, $(X'X)^{-1}X'{\bf y}$, is the best equivariant estimator for ${\bf \boldsymbol\beta}$.
\end{thm}

\begin{proof}
Since the parameter space $\Th$ is transitive under group $\bar{G}$, we choose a special parameter point $\th_0=({\bf \boldsymbol\beta},\Sigma_p)=({\bf 0},I_p)$ to
evaluate the risk function. Thus any equivariant estimator of ${\bf \boldsymbol\beta}$,
$\delta({\bf y})=(X'X)^{-1}X'{\bf y}+X_p^{-1}S({\bf y}){\bf \omega}({\bf z})=X_p^{-1}{\bf \bar{y}}+X_p^{-1}S({\bf y}){\bf \omega}({\bf z})$ has a constant
risk as
follows,
\begin{align*}
	R(\delta)&=EL( \delta,\th_0)=E[L(X_p^{-1}{\bf \bar{y}}+X_p^{-1}S({\bf y}){\bf \omega}({\bf z}),\th_0)|{\bf z}]\\
	&=E\{[X_p^{-1}{\bf \bar{y}}+X_p^{-1}S({\bf y}){\bf \omega({\bf z})}]' X_p'X_p[X_p^{-1}{\bf \bar{y}}+X_p^{-1}S({\bf y}){\bf \omega({\bf z})})]\}\\
	&=E^{\bf z}[E^{{\bf y}|{\bf z}}({\bf \bar{y}}'{\bf \bar{y}})]+2{\bf \omega({\bf z})}'E(S({\bf y})X_p{\bf \bar{y}})+{\bf \omega({\bf z})}'E^{{\bf y}|{\bf
z}}(S^2({\bf y})){\bf \omega({\bf z})}],
\end{align*}
where the minimum is attained at ${\bf \omega}^*={\bf 0}$ as ${\bf \bar{y}},S({\bf y})$ and ${\bf z}$ are pairwise independent. Therefore, the best
estimators is $\delta^*({\bf y})=(X'X)^{-1}X'{\bf y}$.
\end{proof}

The result can be further generalized to any $X$, $n\times p$, $n>p$ matrix with the rank of $p$.
\begin{thm}
\label{MRECoeff2}
For the fixed-$X$ linear model with $X$, an $n\times p$, $n>p$ matrix with the rank of $p$, the least square, $(X'X)^{-1}X'{\bf y}$, is the best equivariant estimator for ${\bf \boldsymbol\beta}$.
\end{thm}

\begin{proof}
Consider the QR decomposition of $X=QX_p$ with $Q$, $n\times p$, $Q'Q=I_p$ and $X_p$, $p\times p$ with the rank of $p$.

Then let ${\bf z}=Q'{\bf y}=X_p\beta+Q'{\bf \epsilon}$ and ${\bf z}$ can be shown to have the form \ref{NormalModel}. Therefore, based on Theorem \ref{MRECoeff}, the MRE estimator is $\delta({\bf z})=(X_p'X_p)^{-1}X_p'{\bf z}=(X'X)^{-1}X'{\bf y}$.
\end{proof}

Notice that here the covariance matrix $\Sigma$ can best be of the form in \ref{NormalModel}, i.e., with at most $p$ distinct values in its diagonal elements. However, such a form in $X$ is against the general logic requirement of equivariance criterion and thus should be treated as a random-$X$ model.

\subsection{\texorpdfstring{Estimation of the Diagonal Covariance Matrix $\Sigma$}{Estimation of the Diagonal Covariance Matrix Sigma}}

Since $\Sigma$ is diagonal with only $p<n$ parameters, it is equivalent to estimate the diagonal matrix, $\Sigma_p$. It is noteworthy that $\Sigma_p$ is
estimable only when $n\geq 2p$. Consider two
widely discussed loss functions: the quadratic loss, $L_q(D,\Sigma_p)=tr((D-\Sigma_p)\Sigma_p^{-1}(D-\Sigma_p) \Sigma_p^{-1})$, and the likelihood loss,
$L_l(D,\Sigma_p)=tr(D\Sigma_p^{-1})-log
|D\Sigma_p^{-1}|-p$, where $D$ is a positive definite $p\times p$ diagonal matrix. The same transformation group $G$ is used in the preservance of the model
with $\bar{g}(\Sigma_p)=C_p\Sigma_pC_p'$ and
both loss functions are invariant under $\tilde{G}_{\Sigma_p}$ with $\tilde{g}(D)=C_p DC_p$ = $C_p^2 D$.
Denote $W=\text{diag}((n_1-1)/(n_1+1),\ldots,(n_p-1)/(n_p+1))$ and one can show the main result as follows.

\begin{thm}
\label{MRECovar}
For the fixed-case linear model in \ref{NormalModel}, $WS^2$ and $S^2$ are the MRE estimators for $\Sigma_p$ under $L_q$ and $L_l$, respectively.
\end{thm}

\begin{proof}
Analog to the univariate case, under loss function $L_q$ and $L_l$, any equivariant estimator $\delta({\bf y})$ of $\Sigma_p$ can be characterized as
\begin{equation}
	\Delta({\bf y})=SH({\bf z})S'=H({\bf z})S^2,
\end{equation}
where $S^2=Diag(s_1^2,\ldots,s_p^2)$, $H({\bf z})=\text{diag}(h_1({\bf z}),\ldots,h_p({\bf z}))$ and ${\bf z}=(\bold{z}_1',\ldots,\bold{z}_p')'$, with $\bold{z}_i=\\((Y_{N_{i-1}+2}-Y_{N_i})/(Y_{N_{i-1}+1}-Y_{N_i}),\ldots,(Y_{N_i-1}-Y_{N_i})/(Y_{N_{i-1}+1}-Y_{N_i}),\text{sigan}(Y_{N_{i-1}+1}-Y_{N_i}))'$. The proof of
equation (2.2) can be found in Appendix \ref{2.2proof}, where it is also shown that $S^2$ and ${\bf z}$ are independent (see Appendix Proofs \ref{s2 and z}).

Since the parameter space $\Th$ is transitive under group $\bar{G}$, we choose a special parameter point $({\bf \boldsymbol\beta},\Sigma_p)=({\bf 0},I_p)=\theta_0$ to
evaluate the risk function.

Firstly, under $L_q$, the constant risk of any equivariant estimator can be calculated as follows,
\begin{align*}
	R(\Delta)&=EL_q(\Delta,\th_0)=E^{\bf z}E^{{\bf y}|{\bf z}}[L_q(H({\bf z})S^2,\th_0)|{\bf z}]\\
	&=E^{\bf z}E^{\bf y}L_q(HS^2,\th_0)=E^{\bf z}\{E_{\th_0}[tr(HS^2-I_p)^2]\}\\
	&=E^{\bf z}\sum_{i=1}^p E_{\th_0}(h_is_i^2-1)^2\\
	&=E^{\bf z}\sum_{i=1}^p (h_i^2E_{\th_0}(s_i^4)-2h_iE_{\th_0}(s_i^2)+1).
\end{align*}
The risk attains the minimum at the point $h_i^*=\frac{E_{\th_0}(s_i^2)}{E_{\th_0}(s_i^4)}=\frac{(n_i-1)^2}{2(n_i-1)+(n_i-1)^2}=\frac{n_i-1}{n_i+1}$ for
$i=1,\ldots,p$ or equivalently,
$H^*=\text{diag}(\frac{n_1-1}{n_1+1},\ldots,\frac{n_p-1}{n_p+1})=W$. Thus, $\delta^*=H^*S^2=WS^2$ is the MRE estimators under $L_q$.

A similar result can be obtained for $L_l$ and the MRE estimators is $\delta^*=S^2$ under $L_l$.
\end{proof}

It is noted that \citet{WuYang2002} has considered such a problem under the common location-scale transformation $g_c$ and the loss function $L_q$ and shown
that no best equivariant estimator exists.

In essence, Theorem \ref{MRECovar}
indicates that the multivariate MRE is a vector of the univariate MREs as in \citet{LehmannCasella1998}, which are sample variances multiplied by a constant. Similar results can be obtained for $\Sigma$ with unknown variances but a known correlation matrix. When the covariance matrix is fully unknown, it is not estimable and there will be no sufficient and complete statistic and thus the above proof won't work.

\renewcommand{\thesection}{\arabic{section}.}
\section{Alternative Proofs via the Principle of Functional Equivariance}
\renewcommand{\thesection}{\arabic{section}}
The above problems can also be solved by converting them into problems of estimation on $p$ independent normal samples if one denote that $\boldsymbol{\mu}=X_p\boldsymbol{\beta}$, as $X_p$ is known and of full rank. In this section, we will discuss
alternative proofs and demonstrate the usage of the principle of functional equivariance.

\begin{thm}
\label{MRECoeff3}
For the invariant problem of estimating $\boldsymbol{\mu}=(\mu_1,\cdots,\mu_p)$ on $p$ independent normal samples $\bold{y}_i \sim N(\mu_i\bold{1}_{n_i}, \sigma_iI_{n_i})$, $i=1,\cdots,p$, under the loss
$L_{\boldsymbol{\mu}}(\bold{d},\boldsymbol{\mu})=(\bold{d}-\boldsymbol{\mu})'\Sigma_p^{-1}(\bold{d}-\boldsymbol{\mu})$ with $\Sigma_p=\text{diag}(\sigma_1,\cdots,\sigma_p)$, and the transformation group $G$, the best equivariant rule is given by ${\bf \bar{y}}=(\bar{Y}_1,\cdots,\bar{Y}_p)'$.
\end{thm}

\begin{proof}
Let $\bold{y}=(\bold{y}_1',\cdots,\bold{y}_p')'$. Then, for any $g\in G$, we have
$$
g(\bold{y})=C\bold{y}+\bold{a},
$$
where
$
C=\text{diag}(c_1I_{n_1},c_2I_{n_2},\cdots,c_pI_{n_p})
$
with $c_i>0$, $i=1,\cdots,p$, and
$
a=(a_1\bold{1}_{n_1}',a_2\bold{1}_{n_2}',\cdots,a_p\bold{1}_{n_p}')'
$
with $a_i\in\mathbb{R}$, $i=1,\cdots,p$. It can be equivalently expressed by that for each $i$,
$$
g_i(\bold{y}_i)=c_i\bold{y}_i+a_i\sim N_{n_i}((c_i\mu_i+a_i)\bold{1}_{n_i},c_i^2\sigma_i^2I_{n_i}),\quad i=1,\cdots,p.
$$
The corresponding transformation group $\bar{G}=\{\bar{g}\}$ on the parameter space can be concluded that
$$
\bar{g}(\boldsymbol{\mu},\Sigma_p)=(C_p\boldsymbol{\mu}+\bold{a}_p,~~C_p\Sigma_pC_p'),
$$
where $C_p=\text{diag}(c_1,\cdots,c_p)$ with $c_i>0$, $i=1,\dots,p$, and $\vec{a}_p=(a_1,\cdots,a_p)'\in\mathbb{R}^p$.

To have the invariance of the loss function, $L_{\boldsymbol{\mu}}(\tilde{g}(\bold{d}),\bar{g}(\boldsymbol{\mu}))=L_2(\bold{d},\boldsymbol{\mu})$, one can show  that
$$
\tilde{g}(\bold{d})=C_p\bold{d}+\bold{a}_p.
$$
Thus, any equivariant estimators should satisfy
$$
\delta_2(g(\bold{y}))=\delta_2(C\bold{y}+\bold{a})=C_p\delta(\bold{y})+\bold{a}_p=\tilde{g}(\delta(\bold{y})).
$$
This implies that $\bar{\bf{y}}$ is equivariant. Furthermore, a characterization of all equivariant estimators is of the form
$$
\delta_2(\bold{y})=\bar{\bold{y}}+S(\bold{y})w(\bold{z}),
$$
where $S(\vec{y})=\text{diag}(s_1,\cdots,s_p)$ with $s_i^2=\frac{1}{n_i-1}\sum_{j=1}^{n_i}(Y_j-\bar{Y}_i)^2$.

Then the risk of an equivariant estimator will be minimized at
\begin{align*}
R_2(\delta_2)=&E_{\boldsymbol{\mu},\Sigma_p}[L_2(\delta_2(\bold{y},\boldsymbol{\mu})]
=E_{\boldsymbol{\mu},\Sigma_p}[L_{\boldsymbol{\beta}}(X_p^{-1}\delta_2(\bold{y}),\boldsymbol{\beta})]\\
=&E_{\bold{0},I_p}(\bar{\bf{y}}'\bar{\bf{y}})=\sum_{i=1}^p\frac{1}{n_i},
\end{align*}
with $w^*=0$.
\end{proof}

In the model \ref{NormalModel}, $k(\boldsymbol{\beta})=X_p\boldsymbol{\beta}$ is a 1-1 and onto mapping and in fact, an element of $\bar{G}$, where $C_p=X_p$, and $a_p=0$. Since, for the converted problems, one can show that the optimal solution is
$A={\bf \bar{y}}$ among all equivariant decision rules, then, for the original problem, the corresponding optimal solution will be $B=k^{-1}(A)=X_p^{-1}{\bf \bar{y}}=(X'X)^{-1}X'{\bf y}$.

\begin{proof}
From the proof for Theorem \ref{MRECoeff2}, one can see that equivariant estimators for problems A and B carry a 1-1 and onto correspondence, $\delta_2(\bold{y})=X_p \delta(\bold{y})$ with $R_2(\delta_2)=R(\delta)$. Therefore, the optimality can be transferred from problem A to B.
\end{proof}

However, the above approach does use the principle of functional equivariance in the following form with a tuned choice of the loss function: For a decision problem $P_A$, there is an optimal solution $A$ and an isomorphic problem $P_B=k(P_A)$, where $k$ is a 1-1
and onto mapping, then $k(A)$ will be an optimal solution for $P_B$.

This naturally holds for equivariance, where Theorem \ref{MRECoeff2} is an example, and implicitly for unbiasedness as one can show that for an optimal problem in the context of unbiasedness, when the optimal solution is unique and an invariant structure exists, it will be equivariant. Such a result can be further extended to the risk-unbiasedness introduced by \citet{Lehmann1951}, which is widely used in the field of machine learning.

Such an approach does have its limitations:
\begin{itemize}
    \item if $k$ is non-linear, then the invariance of the loss function will be lost in the transformation process and thus leading to a new ranking of decision rules; E.g., to estimate the standard deviation $\sigma$ via using the result on estimating the variance $\sigma^2$ in a normal population, where $k(\sigma^2)=\sqrt{\sigma^2}$ but $k(A)$ is no longer optimal.

   \item if the chosen $k$ reduces the dimension and thus is no longer 1-1, then such an
   approach doesn't work. E.g., if one uses $k({\bf \beta})={\bold{1}_p\beta}$ with $g_c \in G_c$ instead of $g \in G$ to estimate $\beta$.

\end{itemize}

\renewcommand{\thesection}{\arabic{section}.}
\section{Future Work and Discussion}
\renewcommand{\thesection}{\arabic{section}}
\setcounter{equation}{0}
This paper has served as an initial effort to explore the usage of equivariance criterion in the field of machine learning, where we are interested at which
method yields the optimal solution  and what properties the optimal solution carries, especially in the context of equivariance criterion. We start from the linear model, the simplest and foundational method in machine
learning.

In this paper, we have established that MRE estimators for the coefficient vector and the condensed covariance matrix is the least square and the vector of
the sample variance within each population, respectively. In addition, we have
demonstrated that in our setting \ref{NormalModel}, the least square estimator is essentially the vector of the sample mean within each population. Such a
finding has further solidified their optimality from the perspective of
multivariate normal distribution theory.

The linear model with a full rank design matrix can be of different forms across literature via the number of populations. The commonly used one is of a
single population, whose characteristic is to assume a common
distribution for all the noises. Naturally, in this setup, one would use a single univariate location-scale transformation to apply equivariance criterion.
In this paper, we relax such an assumption and allow $p$
populations to accommodate the usage of a larger transformation group, which is consisted of multivariate location-scale transformations. From the
perspective of experimental design, such a relaxation is quite intuitive:
the design matrix is chosen carefully before sampling and in a sense, those sampling points are independent and each rank consist of a population. Meanwhile
one can argue that $p$ populations should be the maximum number
allowed in a linear model with fixed-$X$. Even though our results can be extended to the cases where $X$ and $\Sigma$ can be of a more general form, e.g. Theorem \ref{MRECoeff2}, it is against the logical requirement of equivariance criterion for the linear model with a fixed-$X$, where a maximum $p$ populations can be allowed.

In terms of estimating the coefficient vector, the form of $X$ doesn’t
have much impact on the MRE solution and thus \ref{NormalModel} gives a simple path to the MRE solution. However, to estimate the condensed covariance
matrix, it is essential to choose the form and thus decide the size of
each population. In an experimental design setting, one would usually use $n=kp$, where $k$ is an integer, and $n_i=k$. Such a form is recommended both for
its symmetry and computation convenience, which also is a special
case of a seemingly unrelated regression (SUR) problem with a known correlation matrix.

We have discussed the best equivariant estimators for the coefficient vector and the condensed covariance matrix for a normal linear model specially tuned
for equivariance criterion, where the commonly used one is a
special example. Such a model requires a bigger transformation group, which, in turn, will result in a smaller set of equivariant decision rules, where the
MRE estimator exists. Interestingly, \citet{WuYang2002} has shown
that the commonly used single location-scale transformation induces too many equivariant decision rules under $L_q$ that there is no MRE estimator for the
covariance matrix. Meanwhile, each population inside the normal
model is the linear model commonly used in literature and the resulted estimators for each population are the traditional MRE ones, which are equivalent to
the least square solutions.
\citet{KurataMatsuura2016,MatsuuraKurata2020,MatsuuraKurata2021,MatsuuraKurata2024,Matsuura2025Best} used the same transformation group for the SUR model, where a $p$-dimensional distribution family was
considered and samples comes from such a single multivariate population.

The choice of the invariant transformation group is an important topic in equivariance literature. \citet{WuYang2002} presented a case where the group is too
large to allow an optimal solution. Usually, the group is chosen
to be isomorphic to sample space/parameter space, especially when considering the Haar Prior. There are some interesting cases that all invariant
transformation groups pose a nesting relationship and the largest one admits
only the optimal solution for the smaller ones.

Likelihood loss is a multivariate extension of the Stein loss, which is preferred in literature (\citet{Brown1968}). It can be seen that it induces an MRE \&
UMVU estimator, that is always larger than the one under the
quadratic one. Meanwhile, likelihood loss is more evenhanded over the range as the covariance matrix is set to be positive definite. Such a form is quite
similar to the logistic transformation in a generalized linear
model.

An extension to prediction will be a future direction. However, existing frameworks (e.g., \citet{ZhouNayak2015}) on the equivariance criterion can’t handle
the prediction problem well in a linear model. In the literature,
the predicted response is assumed to be unobservable, which is not the case in the linear model. One could also notice that overfitting (\citet{Stone1974})
arises for a prediction problem in a linear model, which usually
doesn't occur in estimation as in deriving the least square solution, sample
prediction error is used, which will converge to a univariate form of $L_{\bf \boldsymbol\beta}$, $({\bf d}-{\bf \boldsymbol\beta})^T X_p^T X_p ({\bf d}-{\bf
\boldsymbol\beta})/\sigma^2$ with $\sigma^2=1$, with an extra term constant to $d$.

Linear model with fixed-$X$ cases though predominantly used in literature, is of limited usage in practice, especially in our settings, where the
experimental design is the ideal scenario. Linear model with random-$X$
cases and mixture cases are more interesting and challenging. The concept of the randomness of linear model has drawn wide attention and numerous efforts
have been spent to clarify the differences between fixedness and
randomness. \citet{Little2019} has given a straightforward definition of randomness as being unknown from a Bayesian view. In his argument, the treatment
indicator from a clinical trial can be considered both fixed and
random. It is true that those semi-controlled covariates pose challenges to the definition of randomness. Individually, it is unknown and thus can be
considered random. Population-wise, its distribution is usually under
control and thus can be treated as a fixed effect in analyses.

The Gauss-Markov Theorem is the fundamental result for the linear model, where the optimality of the least square solution has been established. In most
textbooks, its proof is based on the predominantly assumed fixed-$X$
cases. \citet{Shaffer1991} has shown some interesting results where the Gauss-Markov Theorem no longer holds for some random-$X$ cases. We will investigate
such a phenomenon in the context of equivariance for the
random-$X$ cases.

In terms of randomness, it is noticed that the setup before sampling is crucial and one can classify $X$ before sampling into following categories: known
values, from a known distribution, from a distribution family with
unknown parameters, totally unknown. In a typical experimental design setting, design factors/parameters are of known values, which, in our settings, refer
to fixed-$X$ cases. In a typical clinical trial setting, the
treatment indicator is from a known distribution. In the classical parametric inference setting, we may assume $X$ from a distribution family with unknown
parameters. For the non-parametric setting, X is usually seen as
totally unknown. The latter three scenarios refer to random-$X$ cases, which will be another future topic.

For a linear model with a non-normal distribution, we will refer to extensions of the current results to the generalized linear model, where the challenges
start from the invariant transformations. In a normal linear
model, one can easily find an invariant location-scale transformation group that leave the parameter space transitive, which facilitates the derivation of
the MRE solutions. This may not be the case in a non-normal linear
model.

The principle of Functional Equivariance is widely used in the literature, e.g., \citet{LehmannCasella1998}, in their deduction of the best equivariant estimator in a linear model, uses a canonical form of $X$ and indicated that the result can be extended to the general case of $X$. The key benefit is to convert a problem into a simpler or existing one and then convert the solution back to the original problem.
In this article, we conduct a study of such an approach and find out that equivariance is a sufficient condition for such an approach to hold. It comes to notice that the uniqueness of the optimal solution may also play an important role, which can connects to the Condition A2 of the loss function (\citet{ZhouNayak2014}). Such a condition on the loss function can be termed as "minima-unique".


\medskip
\noindent {\bf Acknowledgment}.
This work was supported by the Fundamental Research Funds for the Central Universities, Sun Yat-sen University (Grant No. 20lgpy145 \& 2021qntd21); the
Science and Technology Program of Guangzhou Project, Fundamental
and Applied Research Project (202102080175); and RDF of Xi'an Jiaotong-Liverpool University (RDF-23-01-073). The authors report there are no competing interests to declare.

The authors would be grateful to the two anonymous referees for their stimulating comments and suggestions leading to a significant improvement of this work.

\newpage
\nocite{*}
\bibliographystyle{stat}
\bibliography{EquivFixedLinearModel}


\begin{appendix}
\section{Technical proof}\label{appA}
{\bf  Identifiability of the Model}
\begin{proof}
For any two parameter points $({\bf \boldsymbol\beta}_1,\Sigma_{1}), ({\bf \boldsymbol\beta}_2,\Sigma_{2})\in\Th$ and any fixed $X$ with full rank, the
density of $N_n(X{\bf \boldsymbol\beta}_1,\Sigma_{1})$ equals to that of $N_n(X{\bf \boldsymbol\beta}_2,\Sigma_{2})$ if and only if $({\bf
\boldsymbol\beta}_1,\Sigma_{p1})=({\bf \boldsymbol\beta}_2,\Sigma_{p2})$, where $\Sigma_i=\left[ \begin{matrix}
\sigma _{i1}^2I_{n_1}&		&		\\
&		\ddots&		\\
&		&		\sigma _{p_i}^2I_{n_{p_i}}\\
\end{matrix} \right] $ for $i=1,2$. Thus, the model is identifiable with respect to $({\bf \boldsymbol\beta},\Sigma_p)$.
\end{proof}

{\bf $G$ being a Group}
\begin{proof}
We aim to prove that
\be
\begin{split}
	G&=\left\{g:g({\bf y})=C{\bf y}+{\bf a}\right\}\\
	\text{with } {\bf a}=(a_1,\ldots,a_p,\ldots,a_p)', C&=Diag(c_1,\ldots,c_p,\ldots,c_p), a_i\in \mathbb{R}\text{ and }c_i>0,i=1,\ldots,p
\end{split}
\enq
satisfies the definition of a group.

(i) {\bf Closedness:} For any two transformations $g_1,g_2\in G$ with
\begin{equation*}
	g_1({\bf y})=C_1{\bf y+a_1},\quad g_2({\bf y})=C_2{\bf y+a_2},
\end{equation*}
where $C_i=Diag(c_{i1},\ldots,c_{ip},\ldots,c_{ip})$, and ${\bf a}_i=(a_{i1},\ldots,a_{ip},\ldots,a_{ip})'$, $i=1,2$,
we have that
\begin{equation*}
	g_2g_1({\bf y})=C_2C_1{\bf y}+C_2{\bf a_1}+{\bf a_2},
\end{equation*}
with $C^*=C_1C_2=Diag(c_{11}c_{21},\ldots, c_{11}c_{21}, \ldots\ldots, c_{1p}c_{2p},\ldots,c_{1p}c_{2p})$ and ${\bf a^*}=C_2{\bf a_1}+{\bf
a_2}=\\(c_{21}a_{11}+a_{21},\ldots, c_{21}a_{11}+a_{21}, \ldots\ldots, c_{2p}a_{1p}+a_{2p},\ldots,c_{2p}a_{1p}+a_{2p})$.
Since $c_i^*=c_{1i}c_{2i}>0$ and $a_i^*=c_{2i}a_{1i}+a_{2i}\in\mathbb{R}$, $i=1,\ldots,p$, we can find that $g_2g_1\in G$.

$(ii)$ {\bf Combination Law:} For any three transformations $g_1,g_2,g_3\in G$ with
\begin{equation*}
	g_i({\bf y})=C_i{\bf y+a}_i,\quad i=1,2,3,
\end{equation*}
we have that
\begin{align*}
(g_1g_2)g_3({\bf y})&=C_1C_2\left(C_3{\bf y}+{\bf a}_3\right)+C_1{\bf a}_2+{\bf a}_1\\
&=C_1(C_2C_3{\bf y}+C_2{\bf a}_3+{\bf a}_2)+{\bf a}_1=g_1(g_2g_3)({\bf y}).
\end{align*}

$(iii)$ {\bf Unit Element:} The transformation $e\in G$ with $C=I_{n},{\bf a=0}$ is the unit element.

$(iv)$ {\bf Inverse Element:} For any transformation $g\in G$, its inverse transformation is $g^{-1}({\bf y})=C^{-1}{\bf y}-C^{-1}{\bf a}$.
\\\\
\end{proof}

{\bf $\bar{G}$ being a Group}
\begin{proof}
We aim to prove that
\begin{equation}
	\begin{split}
		\bar{G}=\{\bar{g}:\bar{g}({\bf \boldsymbol\beta},\Sigma_p)=(X_p^{-1}C_pX_p{\bf \boldsymbol\beta}+X_p^{-1}{\bf 	a}_p,C_p\Sigma_pC_p')\}\\
		\text{with }C_p=Diag(c_1,\ldots,c_p), {\bf a}_p=(a_1,\ldots,a_p)'
	\end{split}	
\end{equation}
satisfies the definition of a group.

$(i)$ {\bf Closedness:} For any two transformations $\bar{g}_1,\bar{g}_2\in\bar{G}$ with
\begin{equation*}
	\bar{g}_i({\bf \boldsymbol\beta},\Sigma_p)=(X_p^{-1}C_{pi}X_p{\bf \boldsymbol\beta}+X_p^{-1}{\bf 	a}_{pi},C_{pi}\Sigma_pC_{pi}'),\quad i=1,2,
\end{equation*}
where $C_{pi}=Diag(c_{i1},\ldots,c_{ip})$, ${\bf a}_{pi}=(a_{i1},\ldots,a_{ip})'$, we have that
\begin{equation*}
	\bar{g}_2\bar{g}_1({\bf \boldsymbol\beta},\Sigma_p)=\left(X^{-1}_pC_{p2}C_{p1}X_p{\bf \boldsymbol\beta}+X_p^{-1}(C_{p2}{\bf a}_{p1}+{\bf
a}_{p2}),C_{p2}C_{p1}\Sigma_p(C_{p2}C_{p1})'\right)
\end{equation*}
with $C_p^*=C_{p2}C_{p1}=Diag(c_{11}c_{21},\ldots,c_{1p}c_{2p})$, ${\bf a}^*_p=C_{p2}{\bf a}_{p1}+{\bf
a}_{p2}=(c_{21}a_{11}+a_{21},\ldots,c_{2p}a_{1p}+a_{2p})$. Since $c_i^*=c_{1i}c_{2i}>0$ and $a_i^*=c_{2i}a_{1i}+a_{2i}\in\mathbb{R}$, $i=1,\ldots,p$, we can
find that $\bar{g}_2\bar{g}_1\in\bar{G}$.

$(ii)$ {\bf Combination Law:} For any three transformations  $\bar{g}_1,\bar{g}_2,\bar{g}_3\in\bar{G}$ with
\begin{equation*}
	\bar{g}_i({\bf \boldsymbol\beta},\Sigma_p)=(X_p^{-1}C_{pi}X_p{\bf \boldsymbol\beta}+X_p^{-1}{\bf 	a}_{pi},C_{pi}\Sigma_pC_{pi}'),\quad i=1,2,3,
\end{equation*}
we have that
\begin{align*}
	&(\bar{g}_1\bar{g}_2)\bar{g}_3({\bf \boldsymbol\beta},\Sigma_p)\\
	=&(X_p^{-1}C_{p1}C_{p2}X_p(X_p^{-1}C_{p3}X_p{\bf \boldsymbol\beta}+X_p^{-1}{\bf a}_{p3})+X_p^{-1}(C_{p1}{\bf a}_{p2}+{\bf
a}_{p1}),(C_{p1}C_{p2})C_{p3}\Sigma_pC_{p3}'(C_{p1}C_{p2})')\\
	=&(X_p^{-1}C_{p1}X_p(X_p^{-1}C_{p2}C_{p3}X_p{\bf \boldsymbol\beta}+X_p^{-1}(C_{p2}{\bf a}_{p3}+{\bf a}_{p2}))+X_p^{-1}{\bf
a}_{p1},C_{p1}(C_{p2}C_{p3})\Sigma_p(C_{p2}C_{p3})'C_{p1}')\\
	=&\bar{g}_1(\bar{g}_2\bar{g}_3)({\bf \boldsymbol\beta},\Sigma_p).
\end{align*}

$(iii)$ {\bf Unit Element:} The Transformation $\bar{e}\in\bar{G}$ with $C_p=I_p$, ${\bf a}_p={\bf 0}$ is the unit element.

$(iv)$ {\bf Inverse Element:} For any transformation $\bar{g}\in\bar{G}$, its inverse transformation is $\bar{g}^{-1}({\bf
\boldsymbol\beta},\Sigma_p)=(X_p^{-1}C_p^{-1}X_p{\bf \boldsymbol\beta}-X_p^{-1}C_p^{-1}{\bf a}_p,C_p^{-1}\Sigma_p(C_p^{-1})')$.
\\\\
\end{proof}

{\bf Transitivity of the Parameter Space: }
For any two parameter points $({\bf \beta_1}, \Sigma_{p1}), ({\bf \beta_2}, \Sigma_{p2}) \in \Theta$, there exists a transformation $\bar{g} \in \bar{G}$
such that $\bar{g}(\beta_1, \Sigma_{p1}) = (\beta_2, \Sigma_{p2})$.

\begin{proof}
Let $g_p({\bf y})=C_p{\bf y}+{\bf a_p}$ with
   \[ C_p = \Sigma_{p2}^{1/2}\Sigma_{p1}^{-1/2} = \text{Diag}(\sqrt{\sigma_{21}^2/\sigma_{11}^2}, \ldots, \sqrt{\sigma_{2p}^2/\sigma_{1p}^2}) \]
and
   \[ {\bf a_p} = X_p{\bf \beta_2} - C_pX_p{\bf \beta_2} .\]

Then
   \[ \bar{g}(\beta_1, \Sigma_{p1}) = (X_p^{-1}\Sigma_{p2}^{1/2}\Sigma_{p1}^{-1/2}X_p\beta_1 + X_p^{-1}(X_p\beta_2 -
   \Sigma_{p2}^{1/2}\Sigma_{p1}^{-1/2}X_p\beta_1), \Sigma_{p2}) \]
   \[ = (\beta_2, \Sigma_{p2}) \]
\end{proof}

{\bf $\tilde{G}_{\bf \boldsymbol\beta}$, $\tilde{G}_{\Sigma_p}$ being Groups}

Analog to the proof for $\bar{G}$,
\begin{equation}
	\begin{split}
		\tilde{G}_{\bf \boldsymbol\beta}=\{\tilde{g}:\tilde{g}({\bf d})&=X_p^{-1}C_pX_p{\bf d}+X_p^{-1}{\bf a}_p\}\text{ and
}\tilde{G}_{\Sigma_p}=\{\tilde{g}:\tilde{g}(D)=C_pDC_p'\}\\
		\text{with }C_p&=Diag(c_1,\ldots,c_p), {\bf a}_p=(a_1,\ldots,a_p)'
	\end{split}
\end{equation}
can be shown to be two groups.
\\\\

{\bf Characterization of Equivariant Estimators for ${\bf \boldsymbol\Sigma_p}$}
\begin{proof}
\label{2.2proof}
It can be easily verified that $\Delta({\bf y})=SH({\bf z})S'=H({\bf z})S^2$ is equivariant under the group action $G_{\Sigma_p}$ if and only if for each of
its diagonal element $\delta_i({\bf y})$, it is equivariant under the transformation group trio $(G,G_c, G_c)$. For $\delta_i({\bf y})$, consider the
transformation $T_{ij}=Y_{N_{i-1}+j}-Y_{N_i}$. Then the problem can be converted into a traditional scale-equivariant one and
the rest follows from Theorem 3.3 in \cite{LehmannCasella1998}.
\end{proof}

{\bf Characterization of Equivariant Estimators for ${\bf \boldsymbol\beta}$}
\begin{proof}
We aim to prove that the estimator $\delta({\bf y})$ is equivariant if and only if it satisfies
\begin{equation*}
	{\bf \delta}({\bf y})=(X'X)^{-1}X' {\bf y}+X_p^{-1}S({\bf y}){\bf \omega}(\bf{z}).
\end{equation*}

To start with, we prove the necessity. From Lemma \ref{LSM}, one can show that the OLS estimator is equivariant. Also, $S(g({\bf
y}))=Diag(c_1s_1,\ldots,c_ps_p)=C_pS({\bf y})$ and ${\bf \omega}({\bf z})$ is invariant under $G$. Thus,
\begin{align*}
	\delta(g({\bf y}))& =X_p^{-1}C_pX_p\cdot (X'X)^{-1}X'{\bf y}+X_p^{-1}{\bf a}_p+X_p^{-1}C_pX_p\cdot X_p^{-1}S({\bf y}){\bf \omega}({\bf z})
    \\
        & =X_p^{-1}C_pX_p \delta({\bf y})+X_p^{-1}{\bf a}_p
     =\tilde{g}(\delta({\bf y})).
\end{align*}
Therefore, $\delta({\bf y})$ is equivariant.

Then we prove the sufficiency. For any equivariant estimator $\delta({\bf y})$, let $\delta_0({\bf y})=\\X_p[\delta({\bf y})-\left(X^{\prime} X\right)^{-1}
X^{\prime}{\bf y}]=X_p[\delta({\bf y})-X_p^{-1}{\bf \bar{y}}]=X_p\delta({\bf y})-{\bf \bar{y}}$ and we have that
\begin{align*}
\delta_0(g({\bf y}))&=X_p\delta(g({\bf y}))- {\bf \overline{g({\bf y})}}\\
&=C_pX_p\delta({\bf y})+{\bf a}_p-C_p{\bf \bar{y}}-{\bf a}_p\\
&=C_p[X_p\delta({\bf y})-{\bf \bar{y}}]=C_p\delta_0({\bf y}).
\end{align*}
Therefore, $\delta_0({\bf y})$ is equivariant under the transformation group trio $(G,G_c, G_c)$.

Similar to the proof above, one can show that $\delta_0({\bf y})$ is equivariant if and only if there exist such an ${\bf \omega}$ that $\delta_0({\bf
y})=S({\bf y}){\bf \omega}(\bf{z})$ as $S({\bf y})$ is equivariant under the transformation group trio $(G,G_c, G_c)$ and $\bf{z}$ is a maximal invariant
under the scale transformation group $G_c$.

Hence, we have ${\bf \delta}({\bf y})=(X'X)^{-1}X' {\bf y}+X_p^{-1}S({\bf y}){\bf \omega}(\bf{z})$.
\end{proof}

\textbf{$(X'X)^{-1}X' {\bf y}$, $S^2$ and ${\bf z}$ are pairwise independent.}
\begin{proof}\label{s2 and z}
Based on Lemma \ref{LSM}, it is easy to show that $(X'X)^{-1}X' {\bf y}$ and $S^2$ are independent.

Next one can show that $({\bf \bar{y}},S^2)$ is complete and sufficient for $(\boldsymbol \eta=X_p \boldsymbol\beta, \Sigma_p)$, ${\bf z}$ is ancillary and
then using Basu Theorem, we will have the independence between $({\bf \bar{y}},S^2)$ and ${\bf z}$.

Alternatively, one can show that all those three as functions based on a linear transformation on ${\bf {y}}$ and the cross-products of their coefficient
matrices are zero matrices.
\end{proof}
\end{appendix}






\end{document}